\documentclass[12pt]{article} 
\usepackage{amssymb}
\usepackage{graphicx}
\usepackage{amsmath}
\usepackage{amsfonts}
\usepackage{hyperref}
\newcommand{\singlespacing}{\let\CS=\@currsize\renewcommand{\baselinestreatch}{1.0}\tiny\CS}
\newcommand{\doublespacing}{\let\CS=\@currsize\renewcommand{\baselinestreatch}{1.5}\tiny\CS}

 \newtheorem{theorem}{Theorem}[section]
 
 \newtheorem{Proposition}{Proposition}[section]

 \numberwithin{equation}{section}
\begin{document}
\begin{center}
\textbf{Characterizations Of Weakly Conformally Flat\\ And Quasi Einstein Manifolds}\\
\vskip.1in
Ramesh Sharma\\
\vskip.1in
University of New Haven, West Haven, CT 06516, USA, E-mail: rsharma@newhaven.edu\\
\end{center}

\begin{abstract} 
\noindent
First, we show that a warped product of a line and a fiber manifold is weakly conformally flat and quasi Einstein if and only if the fiber is Einstein. Next, we characterize and classify contact (in particular, $K$-contact) Riemannian manifold satisfying weakly (and doubly weakly) conformally flat and quasi-Einstein ($\eta$-Einstein) conditions. Finally,  we provide local classification and characterization of a semi-Riemannian (including the 4-dimensional spacetime) with harmonic Weyl tensor and a non-homothetic conformal (including closed) vector field, in terms of Petrov types and Bach tensor.
\end{abstract}
 
\section{Introduction}
In \cite{chen1972hypersurfaces} and \cite{chen1973special}, Chen and Yano introduced the idea of a Riemannian manifold with quasi constant curvature and defined it as a Riemannian manifold ($M,g$) which is conformally flat and whose Ricci tensor $Ric$ satisfies the condition
\begin{equation}\label{1.1}
    Ric = ag + b u\otimes u
\end{equation}
where $a$ and $b$ are smooth functions and $u$ is a 1-form metrically equivalent to a unit vector field $U$ on $M$. Examples of quasi constant curvature manifolds are conformally flat hypersurfaces of a space of constant curvature and canal hypersurfaces (without singularities) of a Euclidean space. We let $g$ be semi-Riemannian and $g(U,U)=\epsilon$, such that $U$ is spacelike if $\epsilon = 1$ and timelike if $\epsilon = -1$. Another example of a quasi constant curvature semi-Riemannian manifold is the warped product $I \times_f F$ of a real open interval $I$ and a Riemannian manifold (fiber) $F$ of constant curvature, with warping function $f>0$ on $I$. In case $I$ is timelike with line-element $-dt^{2}$, $I \times_f F$ is the well known Friedman - Robertson - Walker ($FRW$) cosmological model (see, for example, Hawking and Ellis \cite{hawking1973large}).\\

\noindent
\textbf{Remark:} The equation (1.1) defines a quasi Einstein manifold in the sense of Chaki and Maity \cite{chaki2000quasi}, existence and examples of these manifolds were provided by De and Ghosh \cite{de2004quasi}. In the literature, there is another concept of a quasi Einstein manifold in terms of Bakry-Emery tensor, studied by Case, Shu and Wei \cite{case2011rigidity}.\\

\noindent
In this paper, we follow the former definition [equation (1.1)]. Examples of a quasi Einstein manifold are $\eta$-Einstein Sasakian manifolds (in particular, Sasakian space forms and Heisenberg groups with left invariant metrics as expanding non-gradient Ricci solitons). For details, we refer to Blair \cite{blair2010riemannian}, Chow et al. \cite{Chow} and Ghosh and Sharma \cite{ghosh2014sasakian}. Another example with Lorentzian metric $g$ and a unit timelike vector field $U$, is the perfect fluid solution of Einstein's field equations \cite{hawking1973large}, described by 
\begin{equation*}
    Ric = 4\pi(\rho -  p)g + 8\pi(\rho - 3p)
\end{equation*}
on a 4-dimensional spacetime with energy density $\rho$ and pressure $p$.\\

\noindent
Henceforth we assume the dimension of a semi-Riemannian manifold ($M,g$) to be $> 3$ and denote arbitrary vector fields on $M$ by $X,Y,Z$. Let us weaken the conformal flatness and say that a semi-Riemannian manifold ($M,g$) is weakly conformally flat if and only if
\begin{equation}\label{1.2}
W(X,Y,V) = 0,
\end{equation}
for a non-zero vector field $V$ on $M$, and where $W$ denotes the Weyl tensor of type (1,3).\\

\noindent
In dimension 4, a weakly conformally flat ($M,g$) with a non-null $V$ (i.e. $g(V,V)\ne 0$) is conformally flat. This follows from the following result of Eardley et al. \cite{eardley1986homothetic}. ``If a 4-dimensional semi-Riemannian manifold has its Weyl tensor annihilating a non-null (spacelike or timelike) vector field, then the Weyl tensor vanishes identically". Though this result has been proved for a spacetime, it holds in general, for any semi-Riemannian manifold. The proof is based on symmetry properties, tracelessness and combinatorial computations.\\

\noindent
\textit{In this paper we provide characterizations and classifications of manifolds that are either weakly conformally flat or quasi Einstein, or both within the frameworks of warped product, Lorentzian and contact geometries}. First, we establish the following result.
\begin{theorem}
The semi-Riemannian warped product $(M,g) = I \times_f F$ of dimension $n \geq 4$, with metric $g = \epsilon dt^{2} + f^{2}(t)g_F$ ($g_F$ being the metric on a Riemannian manifold $F$) is weakly conformally flat and quasi-Einstein if any one of the following conditions hold: (i) $g_F$ is Einstein, (ii) $W(x,U,U,y) = 0$ for arbitrary vector fields $x$ and $y$ tangent to $F$ and $U$ the lift of $\partial_t$ ($t$ being the standard coordinate on $I$ (spacelike or timelike)), (iii) Weyl tensor of $g$ is harmonic. Under any one of these conditions, $g$ is Bach flat. 
\end{theorem}
\textbf{Remark:} It follows from the above theorem that an example of a weakly conformally flat and quasi-Einstein is $I \times_f F$ with Einstein fiber $F$. In the 4-dimensional Lorentzian case ($\epsilon = -1$), it is the Friedman - Robertson - Walker ($FRW$) cosmological model, as mentioned earlier.\\

\noindent
Next, we recall from Blair \cite{blair2010riemannian} that a conformally flat contact Riemannian manifold whose Reeb vector field $\xi$ is everywhere an eigen vector of its Ricci tensor, is either of constant curvature 1, or flat (in which case the dimension is 3). Prompted by this result, we consider a $K$-contact manifold $M$ (for which we know that $\xi$ is Killing and is an eigen vector of the Ricci tensor) and check the impact of the doubly weakly conformally flat condition $W(X, \xi, \xi) = 0$ and weakly conformally flat condition $W(X,Y)\xi = 0$ . More precisely, we prove the following result.\\
\begin{Proposition}
A $K$-contact manifold $M$ satisfies the doubly weakly conformally flat condition $W(X, \xi, \xi) = 0$ if and only if $M$ is $\eta$-Einstein. If a $K$-contact manifold satisfies the weakly conformally flat condition $W(X,Y)\xi = 0$, then it is Sasakian and $\eta$-Einstein (i.e. quasi Einstein with $U=\xi$).
\end{Proposition}
Now we provide a classification of a contact Riemannian manifold that is weakly conformally flat and quasi-Einstein with $U$ as the Reeb vector field, by establishing the following result.\\
\begin{theorem}
If a contact Riemannian manifold is weakly conformally flat and quasi-Einstein with $U$ as the Reeb vector field $\xi$, then it is Sasakian, flat (3-dimensional), or locally isometric to $SU(2)$ or $SL(2,\mathbb{R})$ with a left invariant metric.
\end{theorem}

\section{Preliminaries}
Let us briefly review the basic formulas of the $n$-dimensional semi-Riemannian warped product $M = I \times_f F$ of an open real interval $I$ and an ($n-1$)-dimensional semi-Riemannian manifold $F$ (fiber), with warping function $f>0$ on $I$ such that the metric $g$ on $M$ is given by
\begin{equation*}
    g = \pi^{*}(g_I) + (f \circ \pi)^{2} \sigma^{*}(g_F)
\end{equation*}
where $\pi$ and $\sigma$ are projections from  $I \times_f F$ onto $I$ and $F$ respectively, $g_I$ and $g_F$ are the metrics on $I$ and $F$. Locally, we express $g$ as 
\begin{equation*}
    \epsilon dt^{2} + f^{2}(t)(g_F)_{ij}dx^{i}dx^{j}
\end{equation*}
where $\epsilon = 1$ or -1 according as $g(\partial_t, \partial_t) = 1$ or -1, $t$ is the standard coordinate on $I$ and $x^{i}$ are local coordinates on $F$. If $U$ is the lift of $\partial_t$ and $x$, $y$, $z$, $w$ denote arbitrary vector fields tangent to $F$, then we have the following formulas 
\begin{equation}\label{2.1}
    R(x, U)U = -\frac{\ddot{f}}{f}x,  R(x,y)U = 0
\end{equation}
\begin{equation}\label{2.2}
    R(x,U)y = \frac{\epsilon \ddot{f}}{f}g(x,y)U
\end{equation}
\begin{equation}\label{2.3}
g(R(x,y,z),w) = g({^F}R(x,y,z),w) - \epsilon \bigg(\frac{\dot{f}}{f}\bigg)^{2}[g(y,z)g(x,w) - g(x,z)g(y,w)]
\end{equation}
\begin{equation}\label{2.4}
    Ric(U,U) = -\frac{n-1}{f}\ddot{f}
\end{equation}
\begin{equation}\label{2.5}
    Ric(U,x) = 0
\end{equation}
\begin{equation}\label{2.6}
    Ric(x,y) = {^F}Ric(x,y) - \bigg[\frac{\epsilon \ddot{f}}{f} + (n-2)\bigg(\frac{\dot{f}}{f}\bigg)^{2}\bigg]g(x,y)
\end{equation}
\begin{equation}\label{2.7}
    r = \frac{{^F}r}{f^2} - 2\epsilon(n-1)\frac{\ddot{f}}{f} - (n-1)(n-2)\bigg(\frac{\dot{f}}{f}\bigg)^{2}
\end{equation}
where an over-dot means derivative with respect to $t$, and $R$, $Ric$ and $r$ denote the curvature tensor, Ricci tensor and scalar curvature of $M$, while the symbols with a superscript $F$ denote the corresponding quantities of $F$. For details we refer to O'Neill \cite{o1983semi}.\\

\noindent
Secondly, we briefly review contact geometry. A ($2m+1$)-dimensional smooth manifold $M$ is said to be a contact manifold if it carries a smooth 1-form $\eta$ such that $\eta \wedge (d\eta)^m \neq 0$ (where $\wedge$ is the exterior wedge product and $d$ is the exterior derivation). For a given contact form $\eta$, there exists a vector field (Reeb vector field) $\xi$ such that $\eta(\xi)=1$ and $(d\eta)(\xi,.)=0$. Polarizing $d\eta$ on the contact subbundle ($\eta = 0$) we obtain a Riemannian metric $g$ and a (1,1)-tensor field $\varphi$ such that
\begin{equation}\label{2.8}
    (d\eta)(X,Y) = g(X,\varphi Y), \eta(X)=g(\xi, X), \varphi^{2} = -I + \eta \otimes \xi.
\end{equation}
The Riemannian metric $g$ is an associated metric of $\eta$ and ($\varphi, \eta, \xi, g$) is a contact Riemannian structure on $M$. The (1,1)-tensor field $h = \frac{1}{2}\pounds_\xi \varphi$ ($\pounds$ denoting the Lie-derivative operator) on a contact Riemannian manifold is self-adjoint, trace-free and anti-commutes with $\varphi$. The following formulas hold on a contact Riemannian manifold:
\begin{equation}\label{2.9}
    \nabla_ X \xi = -\varphi X - \varphi hX,
\end{equation}
\begin{equation}\label{2.10}
    Ric(\xi,\xi) = 2m - |h|^{2}.
\end{equation}
A contact Riemannian metric is said to be $K$-contact if $\xi$ is Killing, equivalently, $h=0$, or 
\begin{equation}\label{2.11}
    R(X,\xi)\xi = X - \eta(X)\xi
\end{equation}
It follows from (2.10) that $K$-contact metrics are maxima for the Ricci curvature along $\xi$. A contact Riemannian metric on $g$ on $M$ is called Sasakian if the almost Kaehler structure induced on the cone $C(M)$ over $M$ with metric $dr^{2} + r^{2}g$, is integrable (i.e. Kaehler). This Sasakian condition on a contact Riemannian manifold is equivalent to 
\begin{equation}\label{2.12}
    R(X,Y)\xi = \eta(Y)X - \eta(X)Y.
\end{equation}
A Sasakian metric is $K$-contact, however, the converse holds in general, only in dimension 3. A generalization of a Sasakian manifold is the ($k,\mu$)-contact manifold defined as a contact Riemannian manifold $M$ satisfying the nullity condition:
\begin{equation}\label{2.13}
    R(X,Y)\xi = k(\eta(Y)X - \eta(X)Y) + \mu(\eta(Y)hX - \eta(X)hY)
\end{equation}
for real constants $k$, $\mu$. We know that $k \leq 1$, with equality when $M$ is Sasakian. For $k<1$ we have 
\begin{eqnarray}\label{2.14}
    Ric(X,Y) =&& (2m - 2 - m\mu)g(X,Y) + (2m - 2 + \mu)g(hX,Y)\nonumber \\ &+& (m(2k+\mu) - 2m + 2)\eta(X)\eta(Y),
\end{eqnarray}
\begin{equation}\label{2.15}
    r = [2m - 2 + k - m\mu].
\end{equation}
Finally, a contact Riemannian manifold $M$ is said to be $\eta$-Einstein (Boyer, Galicki and Matzeu \cite{boyer2006eta}) if there exist smooth functions $a$ and $b$ on $M$ such that
\begin{equation}\label{2.16}
    Ric = ag + b\eta \otimes \eta.
\end{equation}
For a $K$-contact manifold (in particular, Sasakian manifold) of dimension $>$ 3, we know that $a$ and $b$ are constants. An example of an $\eta$-Einstein manifold is the Sasakian space-form which is a Sasakian manifold of constant $\varphi$-sectional curvature. Sasakian geometry has been extensively studied since its recently perceived relevance in string theory. Sasakian Einstein metrics have received a lot of attention in physics, for example, $p$-brane solutions in superstring theory, Maldacena conjecture (AdS/CFS duality) \cite{Maldacena}.
 
\section{Proofs of The Results}
\textbf{Proof of Theorem 1.1} The Weyl tensor $W$ of type (1,3) is defined by  
\begin{eqnarray}\label{3.1}
    W(X,Y,Z) =&& R(X,Y)Z + \frac{1}{n-2}[Ric(X,Z)Y - Ric(Y,Z)X + g(X,Z)QY\nonumber\\ &-& g(Y,Z)QX]  -\frac{r}{(n-1)(n-2)}[g(X,Z)Y - g(Y,Z)X]
\end{eqnarray}
on an $n$-dimensional semi-Riemannian manifold ($M,g$), where $Q$ is the Ricci operator defined by $g(QX,Y) = Ric(X,Y)$. Using (3.1) and the formulas (2.1) through (2.7) and after a lengthy computation we obtain the following equations
\begin{equation}\label{3.2}
    g(W(x,U,U),y) = -\frac{\epsilon}{n-2} {^F}Ric^{0}(x,y)
\end{equation}
\begin{equation}\label{3.3}
    g(W(x,y,z),U) = 0
\end{equation}
\begin{eqnarray}\label{3.4}
g(W(x,y,z),w) =&& g({^F}W(x,y,z),w) - \frac{1}{(n-2)(n-3)}[g(y,w){^F}Ric^{0}(x,z)\nonumber\\ &-& g(x,w){^F}Ric^{0}(y,z) + g(x,z){^F}Ric^{0}(y,w) \nonumber\\ &-& g(y,z){^F}Ric^{0}(x,w)]
\end{eqnarray}
where $x$, $y$, $z$, $w$ are tangent to $F$, ${^F}W$ the Weyl tensor of $F$ and a superscript $0$ means the traceless part. Parts (i) and (ii) of the theorem follows from (2.6), (3.2) and (3.3). Part (iii) follows from the following result of Gebarowski \cite{gebarowski1992nearly}, "The Weyl tensor of the warped product $I \times_f F$ is harmonic (divergence - free) if and only if $F$ is Einstein". Finally, under any one of the conditions (i), (ii) and (iii), the components of $W$ along $U$ are zero. The Bach tensor $B$ introduced by Bach \cite{Bach} in the context of conformal relativity, is defined as a symmetric second order tensor with components
\begin{equation}\label{3.5}
    B_{ab} = \frac{1}{n-1}\nabla^{c}\nabla^{d}W_{cabd} + \frac{1}{n-2}R^{cd}W_{cabd}
\end{equation}
where $R^{cd}$ are the components of the Ricci tensor. Under any one of the conditions $F$ is Einstein and hence by Gebarowski's result, $W$ is harmonic, also $g$ is quasi-Einstein (implied by $g$ being Weakly conformally flat and quasi Einstein). Using all these facts in (3.5) shows that $B_{ab} = 0$, completing the proof.\\

\noindent
\textbf{Remark:} Taking $x=\partial_i$, $y = \partial_j$ ($\partial_i$ is the coordinate basis $\mathfrak{X}(F)$, we can write $g(W(x,U,U),y)$ as $W_{i00j}$, because $U$ is the lift of $\partial_t$ and $t = x^{0}$ (the coordinate on $I$). The components $W_{i00j}$ are known as the components of the electric part of the Weyl tensor (see \cite{Stephani}) of the spacetime manifold (in which case $\epsilon = -1$).\\

\noindent
\textbf{Proof of Proposition 1.1} Let ($M,\eta,\xi,g$) be a ($2m+1$)-dimensional $K$-contact manifold. Using the definition (3.1) of the Weyl tensor and the $K$-contact property (2.11) we find that
\begin{equation}\label{3.6}
    W(X,\xi)\xi = \frac{r-2m}{2m(2m-1)}[X - \eta(X)\xi] - \frac{1}{2m-1}[QX - 2m\eta(X)\xi]
\end{equation}
If $W(X,\xi)\xi = 0$, then the above equation implies
\begin{equation*}
    QX = \bigg(\frac{r}{2m} - 1\bigg)X + \bigg(2m + 1 - \frac{r}{2m}\bigg)\eta(X)\xi
\end{equation*}
and hence $g$ is $\eta$-Einstein. Conversely, let $g$ be $\eta$-Einstein, i.e. (2.16) holds. Tracing it gives
\begin{equation}\label{3.7}
    r = (2m + 1)a + b
\end{equation}
Writing (2.16) as $QX = aX + b\eta(X)\xi$, substituting $\xi$ for $X$, and using the $K$-contact property: $Q\xi = 2m\xi$ \cite{blair2010riemannian}, we get
\begin{equation}\label{3.8}
    a+b=2m
\end{equation}
Solving (3.7) and (3.8) for $a$ and $b$ and substituting them in (2.16) we get
\begin{equation}\label{3.9}
     QX = \bigg(\frac{r}{2m} - 1\bigg)X + \bigg(2m + 1 - \frac{r}{2m}\bigg)\eta(X)\xi.
\end{equation}
Plugging this value of $QX$ in (3.6) provides $W(X,\xi)\xi  = 0$. This proves the first part. For the second part, we have the stronger hypothesis $W(X,Y)\xi = 0$ on a contact Riemannain manifold $M$. This implies $W(X,\xi,\xi) = 0$, and hence applying the first part, we conclude that $g$ is $\eta$-Einstein, and hence weakly conformally flat and quasi Einstein. Now, a straight forward computation using (3.9) gives
\begin{equation*}
    W(X,Y)\xi = R(X,Y)\xi - \eta(Y)X + \eta(X)Y
\end{equation*}
By our hypothesis, the left hand side is zero, and thus the above equation implies equation (2.12) and hence $g$ is Sasakian. This completes the proof.\\

\noindent
\textbf{Proof of Theorem 1.2} By our hypothesis we have $W(X,Y,\xi) = 0$ and equation (2.16), i.e. $Ric = ag + b\eta \otimes \eta$. Hence $r = (2m + 1)a + b$, by contraction. Using these data we find, after a straight forward computation, that
\begin{equation}\label{3.10}
    R(X,Y)\xi = \frac{a+b}{2m}[\eta(Y)X - \eta(X)Y]
\end{equation}
At this point, we recall the following Schur-type result (Sharma \cite{sharma1995curvature}), "If $R(X,Y)\xi = k(\eta(Y)X - \eta(X)Y)$ for a smooth function $f$ independent of the choice of the vector fields $X$, $Y$ on a contact Riemannian manifold, then $k$ is constant on $M$". Applying this result, in view of equation (3.10), we conclude that $\frac{a+b}{2m}$ is constant, say $k$. Therefore $M$ is a ($k,\mu$)-contact manifold with $\mu = 0$. As $k \leq 1$, we have either $k=1$, in which case $M$ is Sasakian, or $k < 1$. For $k < 1$, comparing (2.16) with (2.14) and noting that $\mu = 0$ and $a+b=2mk$, we have
\begin{equation}\label{3.11}
    (a - 2m + 2)g(X,Y) + (2m - 2 - a)\eta(X)\eta(Y) = 2(m-1)g(hX,Y)
\end{equation}
Contracting it at $X$ and $Y$, and noting that $tr. h = 0$, we find $a=2(m-1)$. Consequently, equation (3.11) reduces to $(m-1)h = 0$. So, $m=1$, because in this case $h\neq0$ (non-Sasakian). Thus, (2.14) reduces to $Ric = 2k\eta \otimes \eta$, i.e. Ricci tensor has rank 1. For $k=0$, $M$ is flat, and for $k\neq0$, $M$ is locally isometric to $SU(2)$ ($k > 0$) or $SL(2, \mathbb{R}$) ($k < 0$) with a left invariant metric (see Blair \cite{blair2010riemannian}). This completes the proof.\\

\section{Concluding Remark}
We may extend the definition (\ref{1.1}) of a quasi-Einstein manifold for a semi-Riemannian metric by allowing the non-zero vector field $U$ to be null (light-like, i.e $g(U,U)=0$).  An example for the null case is the null dust solution of Einstein's field equations, in which the only mass-energy is due to some kind of massless radiation and whose stress-energy tensor is $T=\Phi k \otimes k$ where $k$ is a null vector field specifying the direction of motion of the radiation and $\Phi$ the intensity. This includes plane fronted waves with parallel rays ($pp$-waves) which, in turn include gravitational plane waves. For details, we refer to Stephani et al. \cite{Stephani}.

\end{document}